\documentclass{amsart}
\usepackage{amsmath, amsthm, amssymb, amsfonts}
\usepackage{bm}
\usepackage{dsfont}
\usepackage[utf8]{inputenc}
\usepackage{rotate}
\usepackage{tikz}
\usepackage{tikz-cd}
\usepackage[arrow,matrix,curve]{xy} 	
\usepackage{xcolor}

\theoremstyle{plain}
\newtheorem{theorem}{Theorem}

\theoremstyle{definition}

\newcommand{\norm}[1]{\left\lVert#1\right\rVert}

\newcommand{\RR}{{\mathbb{R}}}

\title[Pair Correlation Statistics and Lattice Point Counting]{Remarks on the Pair Correlation Statistic of Kronecker Sequences and Lattice Point Counting}
\author{Christian Weiss}
\address{{\bf{Ruhr West University of Applied Sciences,}}\\ {{Department of Natural Sciences, Duisburger Str. 100,}}\\{{45479 M\"ulheim an der Ruhr, Germany}}}
\email{christian.weiss@hs-ruhrwest.de}
\date{\today}

\begin{document}

\begin{abstract}
    In this short note, we reformulate the task of calculating the pair correlation statistics of a Kronecker sequence as a lattice point counting problem. This can be done analogously to the lattice based approach which was used to (re-)prove the famous three gap property for Kronecker sequences. We show that recently developed lattice point counting techniques can then be applied to derive that a certain class of Kronecker sequences have $\beta$-pair correlations for all $0 < \beta < 1$. 
\end{abstract}

\maketitle

\section{Introduction}

The concept of Poissonian pair correlations was first introduced in \cite{RS98}. For a sequence $(x_n)_{n \geq 1}$ in the unit interval it means roughly speaking that the number of small distances between points of $(x_n)_{n \geq 1}$ has the expected order of magnitude. Many recent publications, including e.g. \cite{BCC19}, \cite{HKL19}, \cite{LS20}, \cite{TW20}, have addressed the question whether a given sequence (or class of sequences) has Poissonian pair correlations. It seems to be a hard task to find explicit examples of such sequences and only few of them are known, see \cite{BMV15}, \cite{LST21}. Therefore, the notion of $\beta$-pair correlations  for $0 < \beta \leq 1$, where $\beta = 1$ corresponds to Poissonian pair correlations and which is a weaker property than Poissonian pair correlations for $0 < \beta < 1$, recently got attention, see e.g. \cite{LS20}, \cite{WS19}. Most importantly $\beta$-pair correlations were considered in the context of Kronecker sequences $\left\{ n \alpha \right\}_{n \geq 1}$, where $\left\{ \cdot \right\}$ denotes the fractional part of $x \in \mathbb{R}$, and their subsequences.\\[12pt] 
To be more precise, let the norm of a point $x \in \RR$ be defined by $$\norm{x} := \min(x-\lfloor x \rfloor, 1-(x-\lfloor x \rfloor)),$$ where $\lfloor x \rfloor$ denotes the floor bracket. We say that a sequence $(x_n)_{n \geq 1}$ in $[0,1)$ has Poissonian pair correlations if 
$$\lim_{N \to \infty} \frac{1}{N} \# \left\{ 1 \leq l \neq m \leq N \ : \ \norm{x_l - x_m} \leq \frac{s}{N} \right\} = 2s$$
for all $s \geq 0$. Poissonian pair correlations are a generic property of uniformly distributed sequences, i.e. an i.i.d. random sequence (sampled from the uniform distribution in $[0,1)$) almost surely is Poissonian. The opposite is true as well according to \cite{ALP18} and \cite{GL17}. In \cite{NP07}, this notion was generalized to $\beta$-pair correlations for $0 < \beta \leq 1$ if
	$$\lim_{N \to \infty} F_N^\beta(s) := \lim_{N \to \infty} \frac{1}{N^{2-\beta}} \# \left\{ 1 \leq l \neq m \leq N \ : \ \norm{x_l - x_m} \leq \frac{s}{N^\beta} \right\} = 2s$$
holds for all $s \geq 0$. In the literature, $\beta$-pair correlations are also known as number variance, see \cite{Mar07}.\\[12pt]
In this note, we want to contribute to the discussion of the $\beta$-pair correlation property for Kronecker sequences. Our approach uses the lattice based viewpoint of Kronecker sequences which was to our knowledge first introduced in \cite{MS17} and used therein to derive the three gap property of Kronecker sequences, known as Three Gap Theorem. Some of the content of this note might be clear to experts. Nonetheless, our aim is to present the approach to a broader audience. As a direct gain we furthermore obtain from it an easy and short proof that Kronecker sequences $\left\{ n \alpha \right\}_{n \geq 1}$ posses $\beta$-pair correlations for $0 < \beta < 1$ if $\alpha \in \mathbb{R} \setminus \mathbb{Q}$ satisfies some Diophantine condition. Moreover, we hope that the tools presented here will prove to be useful for the analysis of $\beta$-pair correlations also in other contexts, like in higher dimensions, see e.g. \cite{HKL19}, or in an adelic setting, compare \cite{DH21}. Finally, we want to mention that counting the number of lattice points inside a convex set is a classical problem from convex geometry, see e.g. \cite{GW93}. Therefore, also other results from this area (different from the one which we use here), might be useful for calculating the $\beta$-pair correlation statistic.\\[12pt]
In the following, we use the shorthand notation $f \ll_{\gamma} g$ for two functions $f,g$, if there exists a constant $c(\gamma)$ only depending on the parameter $\gamma$ such that $f(x) \leq c(\gamma) g(x)$ holds for all admissible values $x$. If $\lim_{n \to \infty} \frac{f(x)}{g(x)} = 0$ holds we write $f=o(g)$.
\section{Statement of Results}
We now translate the task to calculate the $\beta$-pair correlation statistic $F_N^\beta(s)$ into a lattice point counting problem by amending the ideas in \cite{MS17} to our needs. For that purpose, fix $k \in \mathbb{N}$ and consider the number of points $x_l = \left\{ l\alpha \right\}$ that have distance at most $2s/N^\beta$ from $x_k$. It can be written as
\begin{align*}
\#_{k,N,s,\beta} & = \# \left\{ \norm{(l-k)\alpha + n} \leq s/N^\beta \, | \, (l,n) \in \mathbb{Z}^2, 0 \leq l \leq N \right\}\\ & = \# \left\{ \norm{y - k\alpha} \leq s/N^\beta \vee \, | \, (x,y) \in \mathbb{Z}^2 A_1, 0 \leq x \leq N \right\},
\end{align*}
where
$$A_1 = \begin{pmatrix} 1 & \alpha \\ 0 & 1 \end{pmatrix},$$
i.e. we are counting elements of a certain Bohr sets, compare e.g. to \cite{Cho18}. The expression $\#_{k,N}$ can be further decomposed into
\begin{align*}
1 & + \underbrace{\# \left\{ 0 < y - k\alpha \leq s/N^\beta \, | \, (x,y) \in \mathbb{Z}^2 A_1, 0 < x \leq N \right\}}_{:=\#_{k,N,s,\beta,>}} \\
& + \underbrace{\# \left\{ -s/N^\beta \leq y - k\alpha < 0 \, \vee \, -s/N^\beta \leq y - 1 - k\alpha < 0 \, | \, (x,y) \in \mathbb{Z}^2 A_1, 0 < x \leq N \right\}}_{:=\#_{k,N,s,\beta,<}} ,
\end{align*}
If we think of the Kronecker sequence as a sequence on $\mathbb{R} / \mathbb{Z}$ with ends glued (and equipped with clockwise orientation), the first set may be interpreted as the points succeeding $x_k$  on the unit circle and the second set corresponds to points preceding $x_k$. 
For fixed $s \geq 0$, the $\beta$-pair correlations statistics can thus be formulated as the lattice point counting problem
\begin{align*}
\lim_{N \to \infty} F_N^\beta(s) = \lim_{N \to \infty} \frac{1}{N^{2-\beta}} \sum_{k=1}^N \left( \#_{k,N,s,\beta,>} + \#_{k,N,s,\beta<} + 1 \right),
\end{align*}
if the limit exist. The element $\left\{ l \alpha \right\}$ is predecessor of $\left\{ k \alpha \right\}$ having distance at most $s/N^\beta$ if and only if $\left\{ k \alpha \right\}$ is a successor of $\left\{ l \alpha \right\}$ of distance at most $s/N^\beta$. Hence
\begin{align} \label{eq1}
\lim_{N \to \infty} F_N^\beta(s) = \lim_{N \to \infty} \frac{1}{N^{2-\beta}} \sum_{k=1}^N \left( 2\#_{k,N,s,\beta,>} + 1 \right).
\end{align}
It is well-known for a rather long time that the number of lattice points inside a convex body grows with its volume apart from some error term. However, for classical lattice point counting techniques this error term usually involves the surface of the convex body, see \cite{GW93}. E.g. a prototypical result by Schnell and Wills in \cite{SW91} yields in our situation the estimate
$$V - \frac{2^{3/2}}{2} F \leq \#_{k,N,s,\beta,>} \leq 2V + F + 1,$$
where $V$ denotes the volume of the rectangle and $F$ its surface. The main challenge is thus that the surface of the involved rectangles $F= 2\frac{s}{N^\beta} + 2N$ grows faster in $N$ than the volume $V = sN^{1-\beta}$.\\[12pt] 
While the classical lattice points techniques from convex geometry usually do not impose restrictions on the dimension, more precise error terms can be derived for lower dimensions, where the geometry is more directly accessible in contrast to the general case. In dimension $d=2$, the following result has recently been achieved by Technau and Widmer in \cite{TW20a}. It may be applied to the sets $\#_{k,N,s,\beta,>}$ as long as $N$ is big enough and $\alpha$ satisfies some Diophantine condition (which is the price to be paid for avoiding the surface term). Here, we slightly amend its original formulation such that it fits to our setting.  
\begin{theorem}[Technau, Widmer, \cite{TW20a}, Corollary~1] \label{thm:TW} Let $\phi: (0,\infty) \to (0,1)$ be a non-increasing function such that
$$q|p+q\alpha| \geq \phi(q)$$
for all $(p,q) \in \mathbb{Z} \times \mathbb{N}$. Furthermore assume that $sN^{1-\beta} > 4$ and $0 < s/N^\beta < \sqrt{\alpha}$. Then
$$|\#_{k,N,s,\beta,>} - sN^{1-\beta}| \ll_{\alpha} \frac{\log E}{\phi^{2}(E')},$$
where $E:= \frac{sN^{1-\beta}}{\phi(4N\sqrt{sN^{1-\beta}})}$ and $E' = \sqrt{168}\sqrt{sN^{3-\beta}}E$. 
\end{theorem}
If $\alpha$ is badly approximable, i.e. $\phi(q) = \frac{c}{q}$ for some $c \geq 0$ then we get
$$|\#_{k,N,s,\beta,>} - sN^{1-\beta}| \ll_{\alpha} \log(sN^{1-\beta}).$$
Now assume that $\alpha \in \mathbb{R} \setminus \mathbb{Q}$ is indeed badly approximable and let $s \geq 0$ and $0 < \beta < 1$ be arbitrary. Choose $N$ big enough such that the conditions in Theorem~\ref{thm:TW} are satisfied. For $F_N^\beta(s)$, we hence obtain 
\begin{align*}
\lim_{N \to \infty} F_N^\beta(s) = \lim_{N \to \infty} \frac{1}{N^{2-\beta}} N \cdot \left(2s/N^{1-\beta} + o(\log(sN^{1-\beta})) + 1 \right) = 2s
\end{align*}
This already finishes the proof of our main result which has been proven independently in \cite{WS19} for $\alpha = \frac{1+\sqrt{5}}{2}$ by rather involved number theoretic arguments and in \cite{TW20} by using the $k$-point correlation function.
\begin{theorem} Let $\alpha \in \mathbb{R}$ be a badly approximable number. Then the Kronecker sequence $\left\{ n \alpha \right\}_{n \geq 1}$ has $\beta$- pair correlations for all $0 < \beta < 1$.
\end{theorem}
Note that our proof breaks down for $\beta = 1$ because the first condition in Theorem~\ref{thm:TW} would in this case be $s>4$. Actually, it is easy to see that Kronecker sequences do not have Poissonian by looking at their continued fraction expansion, compare to the simple argument in \cite{WS19}.\\[12pt]
Since the definition of $\beta$-pair correlations depends on the sum over all points $x_k$ only, compare to \eqref{eq1}, it is not directly possible to infer from it results on lattice point counting. However, indeed a stronger property than $\beta$-pair correlations was shown in \cite{WS19} for $\alpha = \frac{1+\sqrt{5}}{2}$, namely
\begin{align} \label{eq2}
\# \left\{ 1 \leq l \leq N \, | \, \norm{x_k-x_l} \leq \frac{s}{N^\beta}  \right\} = 2sN^{1-\beta} + o (2sN^{1-\beta})
\end{align}
for all $0 \leq k \leq N$. This local property can yet again be interpreted as the following result on point counting, namely
\begin{align*}
\left| \# \left\{ (p,q) \right. \right. & \left. \left. \in \mathbb{Z} \times \mathbb{N} \, | \, 0 \leq p + q\alpha - k\alpha \leq 2s/N^{1-\beta}, 0 \leq q \leq N \right\} - 2sN^{1-\beta} \right|\\ & = o (2sN^{1-\beta}).
\end{align*}
Of course, this is a significantly weaker statement than Theorem~\ref{thm:TW} but is still an improvement in comparison to the classical general results from convex geometry.\\[12pt]
\paragraph{Acknowledgment.} I would like to thank Niclas Technau for useful discussions on the topic of this paper.

\bibliographystyle{alpha}
\bibdata{references}
\bibliography{references}

\end{document}